\documentclass[twocolumn]{autart} 

\usepackage[cmex10]{amsmath}
\usepackage{array}
\usepackage{dblfloatfix}
\usepackage{mathtools}
\usepackage{amssymb}
\usepackage{caption}
\usepackage{booktabs}
\usepackage{stackrel}
\usepackage{enumerate}
\usepackage{enumitem}
\usepackage{bbm}
\usepackage[percent]{overpic}
\usepackage{float}
\usepackage{savesym}
\savesymbol{AND}
\usepackage{algorithmic}
\usepackage{algorithm}

\usepackage{soul}
\usepackage[usenames]{color}

\newcommand{\algorithmiconline}{\textbf{Online: }}
\newcommand{\ONLINE}{\STATE \algorithmiconline}
\newcommand{\algorithmicoffline}{\textbf{Offline: }}
\newcommand{\OFFLINE}{\STATE \algorithmicoffline}

\begin{document}

\begin{frontmatter}

\title{Stochastic Model Predictive Control:\\ Output-Feedback, Duality and Guaranteed Performance}

\thanks[footnoteinfo]{Corresponding author R.~R.~Bitmead. The material in this paper was not presented at any conference.}

\author{Martin A. Sehr}\ead{msehr@ucsd.edu},    
\author{Robert R. Bitmead}\ead{rbitmead@ucsd.edu}

\address{Department
    of Mechanical \&\ Aerospace Engineering, University of California, San Diego, La
    Jolla, CA 92093-0411, USA}

\begin{keyword}                        
stochastic control, predictive control, information state, performance analysis, dual optimal control.               
\end{keyword}                           

\begin{abstract}
A new formulation of Stochastic Model Predictive Output Feedback Control is presented and analyzed as a translation of Stochastic Optimal Output Feedback Control into a receding horizon setting. This requires lifting the design into a framework involving propagation of the conditional state density, the \textit{information state}, via the Bayesian Filter and solution of the Stochastic Dynamic Programming Equation for an optimal feedback policy, both stages of which are computationally challenging in the general, nonlinear setup. The upside is that the clearance of three bottleneck aspects of Model Predictive Control is connate to the optimality: output feedback is incorporated naturally; dual regulation and probing of the control signal is inherent; closed-loop performance relative to infinite-horizon optimal control is guaranteed. While the methods are numerically formidable, our aim is to develop an approach to Stochastic Model Predictive Control with guarantees and, from there, to seek a less onerous approximation. To this end, we discuss in particular the class of Partially Observable Markov Decision Processes, to which our results extend seamlessly, and demonstrate applicability with an example in healthcare decision making, where duality and associated optimality in the control signal are required for satisfactory closed-loop behavior.
\end{abstract}

\end{frontmatter}

\section{Introduction}
MPC, in its original formulation, is a full-state feedback law. This underpins two theoretical limitations of MPC: accommodation of output feedback, and extension to include a cogent robustness theory since the state dimension is fixed. This paper addresses the first question. There have been a number of approaches, mostly hinging on replacement of the measured true state by a state estimate, which is computed via Kalman filtering \cite{sehr2016sumptus,yan2005incorporating}, moving-horizon estimator \cite{copp2014nonlinear,sui2008robust}, tube-based minimax estimators \cite{mayne2009robust}, etc. Apart from \cite{copp2014nonlinear}, these designs, often for linear systems, separate the estimator design from the control design. The control problem may be altered to accommodate the state estimation error by methods such as: constraint tightening \cite{yan2005incorporating}, chance/probabilistic constraints \cite{schwarm1999chance}, and so forth.

In this paper, we first consider Stochastic Model Predictive Control (SMPC), formulated as a variant of Stochastic Optimal Output Feedback Control (SOOFC), without regard to computational tractability restrictions. By taking this route, we establish a formulation of SMPC which possesses central features: accommodation of output feedback and duality/probing; examination of the probabilistic requirements of deterministic and probabilistic constraints; guaranteed performance of the SMPC controller applied to the system. Performance bounds are stated in relation to the infinite-horizon-optimally controlled closed-loop performance. We next particularize our performance results to the class of Partially Observable Markov Decision Processes (POMDPs), as is discussed explicitly in~\cite{sehr2017performance}. For this special class of systems, application of our results and verification of the underlying assumptions are computationally tractable, as we demonstrate using a numerical example in healthcare decision making from~\cite{sehr2017tractable}. 

This paper does \textit{not} seek to provide a comprehensive survey of the myriad alternative approaches proposed for Stochastic Model Predictive Control (SMPC). For that, we recommend the numerous available references such as \cite{goodwin2014robust,kouvaritakis2016model,mayne2014model,mesbah2016stochastic}. Rather, we present a new algorithm for SMPC based on SOOFC and prove, particularly, performance properties relative to optimality. As a by-product, we acquire a natural treatment of output feedback via the Bayesian Filter and of the associated controller duality required to balance probing for observability enhancement and regulation. The price we pay for general nonlinear systems is the suspension of disbelief in computational tractability. However, the approach delineates a target controller with assured properties. Approximating this intractable controller by a more computationally amenable variant, as opposed to identifying soluble but indirect problems without guarantees, holds the prospect of approximately attracting the benefits. Such a strategy, using a particle implementation of the Bayesian filter and scenario methods at the cost of losing duality of the control inputs, is discussed in~\cite{sehr2017particle}. Alternatively, as suggested in~\cite{sehr2017tractable}, one may approximate the nonlinear SMPC problem by POMDPs and apply the methods of the current paper directly, resulting in optimality and duality on the approximate POMDP system.
 

\subsection*{Comparison with Other Performance Results}
Our work is related to four central papers discussing performance bounds linking the achieved cost of MPC on the infinite horizon with the cost of infinite-horizon optimal control:
\begin{description}
\item[Gr\"une \&\ Rantzer] 
\cite{grune2008infinite} study the deterministic, full-state feedback situation and provide comparison between the infinite-horizon stochastic optimal cost and the achieved infinite-horizon MPC cost. In particular, the achieved MPC cost is bounded in terms of the computed finite-horizon MPC cost.
\item[Hern\'andes \&\ Lasserre] 
\cite{hernandez1990error} consider the stochastic case with full-state feedback and average as well as discounted costs. Their results yield a comparison between the infinite-horizon stochastic optimal cost and the achieved infinite-horizon MPC cost in terms of the unknown true optimal cost.
\item[Chatterjee \&\ Lygeros] \cite{chatterjee2015stability} also treat the stochastic case with full-state feedback and average cost function. They establish and quantify a bound on the expected long-run average MPC performance 
related to the terminal cost function and its associated monotonicity requirement.
\item[Riggs \&\ Bitmead] \cite{riggs2012mpc} consider the stochastic full-state feedback as an extension to \cite{grune2008infinite} via a discounted infinite-horizon cost function. Similarly to \cite{grune2008infinite}, they establish a performance bound of the achieved infinite-horizon MPC cost in terms of the computed finite-horizon MPC cost.
\item[The current paper] extends \cite{grune2008infinite,riggs2012mpc} to include output feedback stochastic MPC. Achieved performance is bounded in terms of the computed finite-horizon MPC cost. The incorporation of state estimation into the problem is the central contribution.
\end{description}
Each of these works relies on a sequence of assumptions concerning the well-posedness of the underlying optimization problems and specific monotonicity conditions on certain value functions which admit the establishment of stability and performance bounds. 

We summarize the main contribution of this paper, Corollary~\ref{cor:gamma0v2}, for stochastic MPC with state estimation. Subject to cost monotonicity Assumption~\ref{assm:policyw}, which is testable in terms of a known terminal policy
and the terminal cost function
, an upper bound is computable for the achieved infinite-horizon MPC cost in terms of the the computed finite-horizon MPC cost and other parameters of the monotonicity condition. As in \cite{chatterjee2015stability}, we provide an example -- here a POMDP form healthcare -- in which the assumptions are verified, indicating the substance of the assumptions and the nature of the conclusion regarding closed-loop output-feedback stochastic MPC.

\subsection*{Organization of this Paper}
The structure of the paper is as follows. Section~\ref{sec:SOOFC} briefly formulates SOOFC, as used in Section~\ref{sec:SMPC} to present a new SMPC algorithm. After discussing recursive feasibility of this algorithm in Section~\ref{sec:feas}, we proceed by establishing conditions for boundedness of the infinite-horizon discounted cost of the SMPC-controlled nonlinear system in Section~\ref{sec:stab}. Section~\ref{sec:perf} ties the performance of SMPC to the infinite-horizon SOOFC performance. Section~\ref{sec:2} provides a brief encapsulation and post-analysis of the set of technical assumptions in the paper. The results are particularized for POMDPs in Section~\ref{sec:pomdp}, followed by discussion of our numerical example in Section~\ref{sec:example}. We conclude the paper in Section~\ref{sec:disc}. To aid the development, all proofs are relegated to the Appendix.

\subsection*{Notation}
$\mathbb{R}$ and $\mathbb{R}_+$ are real and non-negative real numbers, respectively. The set of non-negative integers is denoted $\mathbb{N}_0$ and the set of positive integers by $\mathbb{N}_1$. We write sequences as $\mathbf{t}^{m} \triangleq \{t_0,t_{1},\ldots,t_{m}\}$, where $m\in\mathbb{N}_0$; $\mathbf{t}^\infty$ is an infinite sequence of the same form. $\operatorname{pdf}(X)$ denotes the probability density function of random variable $X$ while $\operatorname{pdf}(X\vert Y)$ denotes the conditional probability density function of random variable $X$ given jointly distributed random variable $Y$. The acronyms a.s., a.e. and i.i.d. stand for \emph{almost sure, almost everywhere} and \emph{independent and identically distributed}, respectively.
 
\section{Stochastic Optimal Output-Feedback Control}\label{sec:SOOFC}
We consider stochastic optimal control of nonlinear time-invariant dynamics of the form
\begin{align}
x_{k+1}&=f(x_k,u_k,w_k),\quad x_0,\label{eq:state}\\
y_k&=h(x_k,v_k), \label{eq:output}
\end{align}
where $k\in\mathbb{N}_0$, $x_k\in\mathbb{R}^{n_x}$ denotes the state with initial value $x_0$, $u_k\in\mathbb{R}^{n_u}$ the control input, $y_k\in\mathbb{R}^{n_y}$ the measurement output, $w_k\in\mathbb{R}^{n_w}$ the process noise and $v_k\in\mathbb{R}^{n_v}$ the measurement noise. We denote by
\begin{align}\label{eq:init}
\pi_{0\mid -1} &\triangleq \operatorname{pdf}(x_0)
\end{align}
the known a-priori density of the initial state and by
\begin{align*}
\mathbf{\zeta}^k&\triangleq\{y_0,u_0,y_1,u_1,\dots,u_{k-1},y_k\},&\mathbf{\zeta}^0&\triangleq\{y_0\}
\end{align*}
the data available at time $k$. We make the following standing assumptions on the random variables and system dynamics.
\begin{assum}\label{assm:sys}
The dynamics~(\ref{eq:state}-\ref{eq:output}) satisfy
\begin{enumerate}[label=\arabic*.]
\item $f(\cdot,u,\cdot)$ is differentiable a.e. with full rank Jacobian $\forall\,u\in\mathbb{R}^{n_u}$. 
\item $h(\cdot,\cdot)$ is differentiable a.e. with full rank Jacobian.
\item $w_k$ and $v_k$ are i.i.d. sequences with known densities.
\item $x_0, w_k, v_l$ are mutually independent for all $k,l\geq 0$.
\end{enumerate}
\end{assum}
\begin{assum}\label{assm:sysu}
The control input $u_k$ at time instant $k\geq 0$ is a function of the data $\mathbf{\zeta}^k$ and $\pi_{0\mid -1}$.
\end{assum}
As there is no direct feedthrough from $u_k$ to $y_k$, Assumptions~\ref{assm:sys} and~\ref{assm:sysu} assure that system~(\ref{eq:state}-\ref{eq:output}) is a \emph{controlled Markov process}~\cite{BKKUM1986}. Assumption~\ref{assm:sys} further ensures that $f$ and $h$ enjoy the  \textit{Ponomarev 0-property}~\cite{ponomarev1987submersions} and hence that $x_k$ and $y_k$ possess joint and marginal densities.

\subsection{Information State \& Bayesian Filter}
\begin{defn}
The conditional density of state $x$ given data $\mathbf{\zeta}^k$,
\begin{align}\label{eq:pikk}
\pi_{k}&\triangleq\operatorname{pdf}\left(x_{k}\mid \mathbf{\zeta}^k \right),& 
k&\in\mathbb{N}_0,
\end{align}
is the \emph{information state} of system~(\ref{eq:state}-\ref{eq:output}).
\end{defn} 
For a Markov system such as~(\ref{eq:state}-\ref{eq:output}), the information state is propagated via the \emph{Bayesian Filter} (e.g.~\cite{chen2003bayesian,simon2006optimal}):
\begin{align}
\pi_{k}&=
\frac{\operatorname{pdf}(y_{k}\mid x_{k})\,\pi_{k\mid k-1}}{\int \operatorname{pdf}(y_{k}\mid x_{k})\,\pi_{k\mid k-1}\,dx_{k}},\label{eq:BF_rec} \\
\pi_{k+1\mid k} &\triangleq
\int \operatorname{pdf}(x_{k+1} \mid x_{k},u_{k}) \,\pi_{k}\, dx_{k},\label{eq:BF_pred}
\end{align}
for $k\in\mathbb{N}_0$ and density $\pi_{0\mid -1}$ as in~\eqref{eq:init}. For linear dynamics and Gaussian noise, the recursion~(\ref{eq:BF_rec}-\ref{eq:BF_pred}) yields the Kalman Filter.
\begin{defn}\label{def:Tk}
The recursion~(\ref{eq:BF_rec}-\ref{eq:BF_pred}) defines the mapping
\begin{align}\label{eq:Tk}
\pi_{k+1} &= T\left(\pi_{k},y_{k+1},u_k\right),&k&\in\mathbb{N}_0.
\end{align}
\end{defn}

\subsection{Cost and Constraints}
\begin{defn}
$\mathbb{E}_k[\,\cdot\,]$ and $\mathbb{P}_k[\,\cdot\,]$ are expected value and probability with respect to state $x_k$ -- with conditional density $\pi_k$ -- and i.i.d. random variables $\{(w_j,v_{j+1}):j\geq k\}$.
\end{defn}
Given the available data $\mathbf{\zeta}^0$, we aim to select non-anticipatory (i.e. subject to Assumption~\ref{assm:sysu}) control inputs $u_k$ to minimize
\begin{align}\label{eq:costu}
J_N(\pi_{0},\mathbf{u}^{N-1})\triangleq 
\mathbb{E}_0\left[\sum_{j=0}^{N-1}{\alpha^j c(x_j,u_j)} + \alpha^{N}c_{N}(x_{N})\right],
\end{align}
where $N$ is the control horizon, $c:\mathbb{R}^{n_x}\times\mathbb{R}^{n_u}\to\mathbb{R}_+$ the stage cost, $c_N:\mathbb{R}^{n_x}\to\mathbb{R}_+$ the terminal cost and $\alpha\in\mathbb{R}_+$ a discount factor. Drawing from the literature (e.g.~\cite{bertsekas1995dynamic,BKKUM1986}), optimal controls in~\eqref{eq:costu} must inherently be \textit{separated} feedback policies. That is, control input $u_k$ depends on data $\mathbf{\zeta}^k$ and initial density $\pi_{0\mid -1}$ solely through the current information state $\pi_{k}$. Optimality thus requires propagating $\pi_{k}$ and policies $g_k$, where
\begin{align}\label{eq:policies}
u_k = g_k(\pi_{k}).
\end{align}
Cost~\eqref{eq:costu} then reads
\begin{multline}\label{eq:cost}
J_N(\pi_{0},\mathbf{g}^{N-1})= \\
\mathbb{E}_0\left[\sum_{j=0}^{N-1}{\alpha^j c(x_j,g_j(\pi_{j}))} + \alpha^{N}c_{N}(x_{N})\right].
\end{multline}
Extending stochastic optimal control problems with cost~\eqref{eq:cost} to the \emph{infinite horizon} (see~\cite{bertsekas1995dynamic,bertsekas1978stochastic}) typically requires $\alpha < 1$ and omitting the terminal cost term $c_N(\cdot)$, leading to
\begin{align}\label{eq:dcostu}
J_\infty(\pi_{0},\mathbf{g}^{\infty})&\triangleq
\mathbb{E}_0\left[\sum_{j=0}^\infty{\alpha^jc(x_j,g(\pi_{j}))}\right].
\end{align}

In addition to minimizing the expected value cost~\eqref{eq:cost}, we impose probabilistic state constraints of the form
\begin{align}\label{eq:constrx}
\mathbb{P}_k\left[ x_k \in \mathcal{X}_k \right] &\geq 1 - \epsilon_k,&k&\in\mathbb{N}_1
\end{align}
for $\epsilon_k \in[0,1)$. That is, we enforce constraints with respect to the known distributions of the future noise variables and the conditional density of the current state $x_k$, captured by the information state $\pi_k$. Moreover, we consider input constraints of the form
\begin{align}\label{eq:constru}
u_k = g_k(\pi_k) &\in \mathcal{U}_k, &k&\in\mathbb{N}_0.
\end{align}
When discussing infinite-horizon optimal control with cost~\eqref{eq:dcostu}, we replace the state constraints~\eqref{eq:constrx} by the stationary probabilistic state constraints
\begin{align}\label{eq:constrxinf}
\mathbb{P}_k\left[ x_k \in \mathcal{X}_{\infty} \right] &\geq 1 - \epsilon_{\infty},&k&\in\mathbb{N}_1
\end{align}
for $\epsilon_{\infty} \in[0,1)$ and the input constraints~\eqref{eq:constru} by
\begin{align*}
u_k = g_k(\pi_k) &\in \mathcal{U}_{\infty}, &k&\in\mathbb{N}_0.
\end{align*}
\begin{defn}\label{def:Ck}
Denote by $\mathcal{D}$ the set of all densities on $\mathbb R^{n_x}$. Further define $\mathcal C_k\subseteq \mathcal{D}, k\in\mathbb N_1,$ to be the set of all $\pi_k$ of $x_k$ satisfying the probabilistic constraint \eqref{eq:constrx}. Define $\mathcal C_\infty$ likewise for~\eqref{eq:constrxinf}.
\end{defn}

\subsection{Stochastic Optimal Control}
\begin{defn}\label{def:soc}
Given dynamics~(\ref{eq:state}-\ref{eq:output}), $\alpha\in\mathbb{R}_+$ and horizon $N\in\mathbb{N}_1$, define the \emph{finite-horizon stochastic optimal control problem}
\begin{multline*}
\!\!\!\!\mathcal{P}_{N}(\pi_{0}): \hfill
\left\{\!\!\begin{array}{cl}
\inf_{\mathbf{g}^{N-1}}&\!\! J_N(\pi_0,\mathbf{g}^{N-1})\\[0.2cm]
 \text{s.t.}&\!\! \mathbb{P}_j\left[ x_j \in \mathcal{X}_j \right]\geq 1 - \epsilon_j,
\, j=1,\dots,N.\\[0.1cm]
   &\!\! g_j(\pi_{j})\in\mathcal{U}_{j},\, j=0,\dots,N-1.
\end{array}\right.
\end{multline*}
\end{defn}
\begin{defn}\label{def:isoc}
Given dynamics~(\ref{eq:state}-\ref{eq:output}) and $\alpha\in\mathbb{R}_+$, define the \emph{infinite-horizon stochastic optimal control problem}
\begin{align*}
\mathcal{P}_{\infty}(\pi_{0}): 
\left\{ \begin{array}{cl}
\inf_{\mathbf{g}^{\infty}}  & J_{\infty}(\pi_0,\mathbf{g}^{\infty} )\\[0.2cm]
 \text{s.t.}& \mathbb{P}_j\left[ x_j \in \mathcal{X}_{\infty} \right]\geq 1 - \epsilon_{\infty},\, j\in\mathbb{N}_1.\\[0.1cm]
  & g_j(\pi_{j})\in\mathcal{U}_{\infty},\, j\in\mathbb{N}_0.
\end{array}\right.
\end{align*}
\end{defn}
\begin{defn}\label{def:feas}
$\pi_{0}$ is feasible for $\mathcal{P}_{N}(\cdot)$ if there exists a sequence of policies $\mathbf{g}^{N-1}$ such that, $\{w_{j},v_{j+1}\}_{j\geq 0}$-a.s., $u_{j} = g_j(\pi_{j})$ satisfy the constraints and $J_N(\pi_0,\mathbf{g}^{N-1})$ is finite. Define feasibility likewise for $\mathcal{P}_{\infty}(\pi_{0})$.
\end{defn}
In Stochastic Optimal Control, feasibility entails the existence of policies $g_k(\cdot)$ such that for any $\pi_k\in\mathcal C_k$, $g_k(\pi_k)\in\mathcal{U}_k$ and
\begin{align*}\label{eq:constrpi}
\pi_{k+1}=T(\pi_k,y_{k+1},g_k(\pi_k))&\in\mathcal C_{k+1},&&(w_k,v_{k+1})-\text{a.s.}
\end{align*}
Even though the state constraints~\eqref{eq:constrx} are probabilistic, this condition results in an equivalent almost sure constraint on the conditional state densities. The stochastic optimal feedback policies in $\mathcal{P}_N(\pi_0)$ may now be computed in principle by solving the Stochastic Dynamic Programming Equation (SDPE),
\begin{equation}
\begin{aligned}
\label{eq:DP1}
V_k(\pi_{k}) \triangleq 
\inf_{g_k(\cdot)}&\ \mathbb{E}_k \left[ c(x_k,g_k(\pi_k)) +\alpha V_{k+1}(\pi_{k+1})\right], \\
\text{s.t.}\,& \ \pi_{k+1}\in\mathcal{C}_{k+1},\ (w_k,v_{k+1})-\text{a.s.}\\
& \ g_k(\pi_k) \in \mathcal{U}_k
\end{aligned}
\end{equation}
for $k = 0,\ldots,N-1$ and $\pi_k\in\mathcal{C}_k$. The equation is solved backwards in time, from its terminal value
\begin{align}\label{eq:DP2}
V_N(\pi_{N}) &\triangleq \mathbb{E}_N \left[ c_N(x_{N})\right], &
\pi_N \in \mathcal{C}_{N}.
\end{align}

Solution of the SDPE is the primary source of the restrictive computational demands in Stochastic Optimal Control. The reason for this difficulty lies in the dependence of the future information state in each step of~(\ref{eq:DP1}-\ref{eq:DP2}) on the current and future control inputs. While the dependence on future control inputs is limiting even in deterministic control, the computational burden is drastically worsened in the stochastic case because of the complexity of the operator $T_k$ in \eqref{eq:Tk}. On the other hand, optimality via the SDPE leads to a control law of \emph{dual} nature. \textit{Dual optimal control} connotes the compromise in optimal control between the control signal's function to reveal the state and its function to regulate that state. These dual actions are typically antagonistic~\cite{Feldbaum1965}. The duality of stochastic optimal control is a generic feature, although there exist some problems -- called \textit{neutral} -- where the probing nature of the control evanesces, linear Gaussian control being one such case.

Notice that, while the Bayesian Filter~(\ref{eq:BF_rec}-\ref{eq:BF_pred}) can be approximated to arbitrary accuracy using a Particle Filter~\cite{simon2006optimal}, the SDPE cannot be easily simplified without loss of optimal probing in the control inputs. While control laws generated without solution of the SDPE can be modified artificially to include certain excitation properties, as discussed for instance in~\cite{genceli1996new,marafioti2014persistently}, such approaches are suboptimal and do not generally enjoy the theoretical guarantees discussed below. For the stochastic optimal control problems considered here, excitation of the control signal is incorporated automatically and as necessary through the optimization. The optimal control policies, $g^\star_j(\cdot)$, will inherently inject excitation into the control signal depending on the quality of state knowledge embodied in $\pi_k$. 

\section{Stochastic Model Predictive Control}\label{sec:SMPC}
\begin{algorithm}[H]                   
\caption*{(Dual Optimal) SMPC}\label{RHSOC}                       
\begin{algorithmic}[1]
\REQUIRE $\pi_{0\mid -1}\in\mathcal{D}$ and $\alpha\in\mathbb{R}_+$
\OFFLINE Solve $\mathcal{P}_N(\cdot)$ for $g_0^{\star}(\cdot)$ via~(\ref{eq:DP1}-\ref{eq:DP2})
\ONLINE
\FOR{$k\in\mathbb{N}_0$}
\STATE Measure $y_k$
\STATE Compute $\pi_k$
\STATE Apply first optimal control policy, $u_{k} = g_0^{\star}(\pi_{k})$
\STATE Compute $\pi_{k+1\mid k}$
\ENDFOR
\end{algorithmic}
\end{algorithm}
Notice how this algorithm differs from common practice in SMPC \cite{kouvaritakis2015stochastic,mesbah2016stochastic} in that we explicitly use the information states $\pi_k\in\mathcal{D}$. Throughout the literature, these information states -- conditional densities -- are replaced by best available, or certainty-equivalent state estimates in $\mathbb{R}^{n_x}$. While this makes the problem more tractable, one no longer solves the underlying stochastic optimal control problem. As we shall demonstrate in this paper, using information state $\pi_k$ and optimal policy $g_0^{\star}(\cdot)$ resulting from solution of Problem $\mathcal{P}_N(\pi_k)$ at each time instance leads to a number of results regarding closed-loop performance on the infinite horizon. 

\section{Recursive Feasibility}\label{sec:feas}
\begin{assum}\label{assm:0}
$\pi_{0\mid -1}$ yields $\pi_0$ feasible for $\mathcal{P}_N(\cdot)$, $v_0$-a.s.
\end{assum}
\begin{assum}\label{assm:1}
The constraints in $\mathcal{P}_N(\cdot)$ and $\mathcal{P}_{\infty}(\cdot)$, for $j=1,\ldots,N-1$, satisfy
\begin{align*}
\mathcal{C}_{j+1} \subseteq \mathcal{C}_j &\subseteq \mathcal{C}_{\infty},& 
\mathcal{U}_{j} \subseteq \mathcal{U}_{j-1}&\subseteq \mathcal{U}_{\infty}.
\end{align*}
\end{assum}
\begin{assum}\label{assm:2}
For all densities $\pi_k\in\mathcal{C}_N$, there exists a policy $\tilde{g}(\pi_k)$ satisfying
\begin{align*}
\tilde{g}(\pi_k)&\in\mathcal{U}_{N-1},\\
T(\pi_k,y_{k+1},\tilde{g}(\pi_k))&\in\mathcal{C}_N,\ (w_k,v_{k+1})-\text{a.s.},\\
c(x_k,\tilde{g}(\pi_k))&< \infty.
\end{align*}
\end{assum}
\begin{thm}
  \label{thm:srf}
Given Assumptions~\ref{assm:0}-\ref{assm:2}, SMPC yields $\pi_k$ feasible for $\mathcal{P}_N(\cdot)$, $\{w_{j},v_{j+1}\}_{j\geq 0}$-a.s., for all $k\in\mathbb{N}_1$.
\end{thm}
The proof of this result follows directly as a stochastic version of the corresponding result in deterministic MPC, e.g. \cite{grune2011nonlinear}. Notice that recursive feasibility and compact $\mathcal X_1$ immediately implies a stability result independent of the cost~\eqref{eq:cost}, i.e.
\begin{align}\label{eq:cor:stab}
\mathbb{P}_k[x_k\in\mathcal{X}_1] &\geq 1 - \epsilon_1,& \{w_{j},v_{j+1}\}_{j\geq 0}-\text{a.s.},
\end{align}
for $k\in\mathbb{N}_1$.

\section{Convergence and Stability}\label{sec:stab}
\begin{assum}\label{assm:policy}
For a given $\alpha\in\mathbb{R}_+$, the terminal feedback policy $\tilde{g}(\pi)$ specified in Assumption~\ref{assm:2} satisfies
\begin{align}\label{eq:ineq}
\alpha\,\mathbb{E}_{\pi}\left[ c_N(f(x,\tilde{g}(\pi),w)) \right] - c_N(x) \stackrel{\text{a.s.}}{\leq}
 - c(x,\tilde{g}(\pi))
\end{align}
for all densities $\pi$ of $x$ with $\pi\in \mathcal{C}_N$. The expectation $\mathbb{E}_{\pi}[\cdot]$ is with respect to state $x$ -- with density $\pi$ -- and $w$.
\end{assum}
For $\alpha \geq 1$, Assumption~\ref{assm:policy} can be interpreted as the existence of a stochastic Lyapunov function on the terminal set of densities, $\mathcal{C}_N$. If~\eqref{eq:ineq} holds for $\alpha \geq 1$, it naturally holds for all $\alpha \in (0,1]$.
\begin{thm}\label{thm:1}
Given Assumptions~\ref{assm:0}-\ref{assm:policy}, SMPC yields
\begin{align}\label{eq:summability}
\lim_{M\to\infty}\sum_{k = 0}^{M} { \alpha^k c(x_k,g_0^{\star}(\pi_{k})) } \stackrel{\text{a.s.}}{<} \infty.
\end{align}
\end{thm}
While the discount factor $\alpha$ may not seem to play a major role in this result, notice that small values of $\alpha$ may be required to satisfy Assumption~\ref{assm:policy}. For $\alpha\geq 1$, \eqref{eq:summability} implies almost sure convergence to 0 of the achieved stage cost.
\begin{assum}\label{assm:detect}
State $x$ is detectable via the stage cost:
\begin{align*}
c(x_k,u_k)&\stackrel{\text{a.s.}}{\to} 0 \text{ as } k\to\infty&\implies&& 
x_k&\stackrel{\text{a.s.}}{\to}\mathcal{X} \text{ as } k\to\infty.
\end{align*}
\end{assum}
\begin{thm}\label{thm:2}
Given Assumptions~\ref{assm:0}-\ref{assm:detect}, SMPC with $\alpha \geq 1$ yields
\begin{align*}
\lim_{M\to\infty}\sum_{k = 0}^{M} { c(x_k,g_0^{\star}(\pi_{k})) } \stackrel{\text{a.s.}}{<} \infty
\end{align*}
and
\begin{align}\label{eq:converge}
x_k\stackrel{\text{a.s.}}{\to}\mathcal{X},\text{ as } k\to\infty.
\end{align} 
\end{thm}
While~\eqref{eq:converge} holds only for $\alpha \geq1$, notice that SMPC for $\alpha\in[0,1)$ with recursive feasibility possesses the default stability property~\eqref{eq:cor:stab}. For zero terminal cost $c_N(x) \equiv 0$, Assumption~\ref{assm:5} replaces Assumption~\ref{assm:policy} to guarantee~\eqref{eq:summability}, a finite discounted infinite-horizon SMPC cost.
\begin{assum}\label{assm:5}
The terminal feedback policy $\tilde{g}(\pi)$ specified in Assumption~\ref{assm:2} satisfies
\begin{align*}
c(x,\tilde{g}(\pi)) \stackrel{\text{a.s.}}{=} 0
\end{align*}
for all densities $\pi$ of $x$ with $\pi\in \mathcal{C}_N$.
\end{assum}
\begin{cor}\label{cor:0term}
Given Assumptions~\ref{assm:0}-\ref{assm:2} and~\ref{assm:5}, SMPC with zero terminal cost $c_N(x) \equiv 0$ yields
\begin{align*}
\lim_{M\to\infty}\sum_{k = 0}^{M} { \alpha^k c(x_k,g_0^{\star}(\pi_{k})) } \stackrel{\text{a.s.}}{<} \infty.
\end{align*}
Moreover, if $\alpha = 1$ and Assumption~\ref{assm:detect} is added, we have
\begin{align*}
x_k\stackrel{\text{a.s.}}{\to}\mathcal{X},\text{ as } k\to\infty.
\end{align*}
\end{cor}

\section{Infinite-Horizon Performance Bounds}\label{sec:perf}
In the following, we establish performance bounds for SMPC, implemented on the infinite horizon as a proxy to solving the infinite-horizon stochastic optimal control problem $\mathcal{P}_{\infty}(\pi)$. These bounds are in the spirit of previously established bounds reported for deterministic MPC in~\cite{grune2008infinite} and the stochastic full state-feedback case in~\cite{riggs2012mpc}.
\begin{assum}\label{assm:gamma}
There exist $\gamma\in [0,1]$ and $\eta\in\mathbb{R}_+$ such that
\begin{multline}\label{eq:gamma}
\mathbb{E}_0\left[ V_{0}(T(\pi_0,y_1,g_0^{\star}(\pi_0))) - V_{1}(T(\pi_0,y_1,g_0^{\star}(\pi_0))) \right] \leq \\
\gamma\,\mathbb{E}_0\left[ c(x_0,g_0^{\star}(\pi_0)) \right] + \eta
\end{multline}
for all densities $\pi_0$ of $x_0$ which are feasible in $\mathcal{P}_N(\cdot)$. 
\end{assum}
\begin{defn}\label{def:opt}
Denote by $\mathbf{g}^{MPC}$ the SMPC implementation of policy $g_0^{\star}(\cdot)$ on the infinite horizon, i.e. 
\begin{align*}
\mathbf{g}^{MPC} \triangleq \{ g_0^{\star},g_0^{\star},g_0^{\star},\ldots\}.
\end{align*}
Similarly, $\mathbf{g}^{\star^{N-1}}$ and $\mathbf{g}^{\star^\infty}$ are the optimal sequences of policies in Problems $\mathcal{P}_N(\cdot)$ and $\mathcal{P}_\infty(\cdot)$, respectively.
\end{defn}
\begin{thm}\label{thm:bounds}
Given Assumptions~\ref{assm:0}-\ref{assm:2} and~\ref{assm:gamma}, SMPC with $\alpha\in[0,1)$ yields
\begin{multline}\label{eq:bounds}
(1- \alpha\gamma)\, J_{\infty}(\pi_0,\mathbf{g}^{\star^\infty}) \leq \\
(1- \alpha\gamma)\, J_{\infty}(\pi_0,\mathbf{g}^{MPC}) \leq \\
J_{N}(\pi_0,\mathbf{g}^{\star^{N-1}}) + \frac{\alpha}{1 - \alpha}\eta.
\end{multline}
\end{thm}
In the special case $\gamma = 0$, we impose the following assumption on the terminal cost to obtain an insightful corollary to Theorem~\ref{thm:bounds}.
\begin{assum}\label{assm:policyw}
For $\alpha\in[0,1)$, there exists $\eta\in\mathbb{R}_+$ such that the terminal policy $\tilde{g}(\cdot)$ specified in Assumption~\ref{assm:2} satisfies
\begin{multline*}
\mathbb{E}_{\pi}\left[\alpha\, c_N(f(x,\tilde{g}(\pi),w)) - c_N(x) \right] \stackrel{}{\leq}\\
-\mathbb{E}_{\pi}\left[c(x,\tilde{g}(\pi))\right]
 + \frac{\eta}{\alpha^{N-1}},
\end{multline*}
for all densities $\pi$ of $x$ with $\pi\in \mathcal{C}_{N}$. The expectation $\mathbb{E}_{\pi}[\cdot]$ is with respect to state $x$ -- with density $\pi$ -- and $w$.
\end{assum}
\begin{cor}\label{cor:gamma0v2}
Given Assumptions~\ref{assm:0}-\ref{assm:2} and~\ref{assm:policyw}, SMPC with $\alpha\in[0,1)$ yields 
\begin{multline*}
J_{\infty}(\pi_0,\mathbf{g}^{\star^\infty}) \leq 
J_{\infty}(\pi_0,\mathbf{g}^{MPC}) \leq \\
J_{N}(\pi_0,\mathbf{g}^{\star^{N-1}}) + \frac{\alpha}{1 - \alpha}\eta.
\end{multline*}
\end{cor}
This Corollary relates the following quantities: \textit{design cost}, $J_N(\pi_0,\mathbf{g}^{*^{N-1}}),$ which is known as part of the SMPC calculation, \textit{optimal cost} $J_\infty(\pi_0,\mathbf{g}^{*^\infty})$ which is unknown (otherwise we would use $\mathbf{g}^{*^\infty}$), and unknown infinite-horizon SMPC \textit{achieved cost} $J_\infty(\pi_0,\mathbf{g}^{MPC})$.

\section{Analysis of Assumptions}\label{sec:2}
The sequence of assumptions becomes more inscrutable as our study progresses. However, they deviate only slightly from standard assumptions in MPC, suitably tweaked for stochastic applications. Assumptions~\ref{assm:sys} and \ref{assm:sysu} are regularity conditions permitting the development of the Bayesian filter via densities and restricting the controls to causal policies. Assumptions~\ref{assm:0} and \ref{assm:1} limit the constraint sets and initial state density to admit treatment of recursive feasibility.

Assumptions~\ref{assm:2}, \ref{assm:policy}, \ref{assm:5} and \ref{assm:policyw} each concerns a putative terminal control policy, $\tilde{g}(\cdot)$. Assumption~\ref{assm:2} implies positive invariance of the terminal constraint set under $\tilde{g}$. Using the martingale analysis of the proof of Theorem~\ref{thm:2}, Assumption~\ref{assm:policy} ensures that the extant $\tilde{g}$ achieves finite cost-to-go on the terminal set. The cost-detectability Assumption~\ref{assm:detect} is familiar in Optimal Control to make the implication that finite cost forces state convergence. Assumption~\ref{assm:5} temporarily replaces Assumption~\ref{assm:policy} only to consider the zero terminal cost case. Assumptions~\ref{assm:gamma} and \ref{assm:policyw} presume monotonicity of the finite-horizon cost with increasing horizon, firstly for the optimal policy $g_0^{\star}$ and then for the putative terminal policy, $\tilde{g}$ on the terminal set. These monotonicity assumptions mirror those of, for example, \cite{grune2008infinite} for deterministic MPC and \cite{riggs2012mpc} for full-state stochastic MPC. They underpin the deterministic Lyapunov analysis and the stochastic Martingale analysis based on the cost-to-go. These assumptions are validated for a POMDP example in Section~\ref{sec:example}.


\section{Dual Optimal Stochastic MPC for POMDPs}
\label{sec:pomdp}
We now proceed by particularizing the performance results from Section~\ref{sec:perf} for the special class of POMDPs, as suggested for instance in~\cite{sehr2017performance,sehr2017tractable,sunberg2013information}. This class of problems is characterized by probabilistic dynamics on a finite state space $X = \{1,\ldots,n\}$, finite action space $U = \{1,\ldots,m\}$, and finite observation space $Y = \{1,\ldots,o\}$. POMDP dynamics are defined by the conditional state transition and observation probabilities
\begin{align}\label{eq:pomdpx}
\mathbb{P}\left(x_{t+1} = j \mid x_t = i, u_t = a\right) &= p_{ij}^{a}, \\
\label{eq:pomdpy}
\mathbb{P}\left(y_{t+1} = \theta \mid x_{t+1} = j, u_t = a\right) &= r_{j\theta}^{a},
\end{align}
where $t\in\mathbb{N}_0$, $i,j\in X$, $a\in U$, $\theta\in Y$. The state transition dynamics~\eqref{eq:pomdpx} correspond to a conventional Markov Decision Process (MDP, e.g.~\cite{puterman2014markov}). However, the control actions $u_t$ are to chosen based on the known initial state distribution $\pi_0 = \operatorname{pdf}(x_0)$ and the sequences of observations, $\{y_1,\ldots,y_t\}$, and controls $\{u_0,\ldots,u_{t-1}\}$, respectively. That is, we are choosing our control actions in a Hidden Markov Model (HMM, e.g.~\cite{elliott2008hidden}) setup. Notice that, while POMDPs conventionally do not have an initial observation $y_0$ in~\eqref{eq:pomdpy}, as is commonly assumed in nonlinear system models of the form~(\ref{eq:state}-\ref{eq:output}), one can easily modify this basic setup without altering the following discussion.

Given control action $u_t = a$ and measured output $y_{t+1} = \theta$, the information state $\pi_t$ in a POMDP is updated via
\begin{align*}
\pi_{t+1,j} = \frac{\sum_{i\in X} 
\pi_{t,j} p_{ij}^a r_{j\theta}^a}{\sum_{i,j\in X} \pi_{t,j} p_{ij}^a r_{j\theta}^a} ,
\end{align*}
where $\pi_{t,j}$ denotes the $j^{\text{th}}$ entry of the row vector $\pi_t$. To specify the cost functionals~\eqref{eq:cost} and~\eqref{eq:dcostu} in the POMDP setup, we write the stage cost as $c(x_t,u_t) = c_{i}^{a}$ if $x_t = i\in X$ and $u_t = a\in U$, summarized in the column vectors $c(a)$ of the same dimension as row vectors $\pi_k$. Similarly, the terminal cost terms are $c_N(x_t) = c_{i,N}$ if $x_N = i\in X$, summarized in the column vector $c_N$. The infinite horizon cost functional defined in Section~\ref{sec:SOOFC} then follows as 
\begin{align*}
J_\infty(\pi_{0},g) & =
\mathbb{E}_0 \left[ \sum_{k=0}^\infty{\alpha^k \pi_k c(g_k(\pi_k)) } \right],
\end{align*}
with corresponding finite-horizon variant
\begin{align*}
J_N(\pi_{0},g) &= 
\mathbb{E}_0 \left[ \sum_{k=0}^{N-1}{\alpha^k \pi_k c(g_k(\pi_k))} + 
\alpha^{N}\pi_N c_N \right].
\end{align*}
Extending~(\ref{eq:DP1}-\ref{eq:DP2}), optimal control decisions may then be computed via
\begin{multline}\label{eq:DP1pomdp}
J_{N-k}^{\star}(\pi_k) = \min_{g_k(\cdot)} \Bigg\{ \pi_{k} c(g_k(\pi_k)) \\ +
\alpha\sum_{\theta \in {Y}} \mathbb{P}\left(y_{k+1} = \theta \mid \pi_k,\,g_k(\pi_k) \right) 
J_{N-k-1}^{\star}(\pi_{k+1}) \Bigg\}, 
\end{multline}
for $k = 0,\ldots,N-1$, from terminal value function
\begin{align}\label{eq:DP2pomdp}
J_{0}^{\star}(\pi_N) = \pi_N c_N .
\end{align}
\begin{assum}\label{assm:pomdp}
For $\alpha\in[0,1)$, there exist $\eta\in\mathbb{R}_+$ and a policy $\tilde{g}(\cdot)$ such that
\begin{align}\label{eq:pomdpassm}
\mathbb{E}_{0}\left[\alpha\, \pi_1 c_N \right] \stackrel{}{\leq}
\mathbb{E}_{0}\left[ \pi_0 c_N - \pi_0 c(\tilde{g}(\pi_0))  \right]
 + \frac{\eta}{\alpha^{N-1}},
\end{align}
for all densities $\pi_0$ of $x_0\in X$.
\end{assum}
\begin{thm}[\cite{sehr2017performance}]
\label{thm:pomdp}
Given Assumption~\ref{assm:pomdp}, SMPC for POMDPs with $\alpha\in[0,1)$ yields 
\begin{multline*}
J_{\infty}(\pi,\mathbf{g}^{\star^\infty}) \leq 
J_{\infty}(\pi,\mathbf{g}^{MPC}) \leq \\
J_{N}(\pi,\mathbf{g}^{\star^{N-1}}) + \frac{\alpha}{1 - \alpha}\eta,
\end{multline*}
for all densities $\pi$ of $x\in X$.
\end{thm} 
A special case of Corollary~\ref{cor:gamma0v2}, this result allows us to bound the achieved infinite-horizon cost of SMPC on POMDPs. In this special case, we can compute the dual optimal control policies and verify Assumption~\ref{assm:pomdp} numerically, as is demonstrated for a particular example below.

\section{An Example in Healthcare Decision Making}
\label{sec:example}
\subsection{Problem Setup}
Consider a patient treated for a specific disease which can be managed but not cured. For simplicity, we assume that the patient does not die under treatment. While this transition would have to be added in practice, it results in a time-varying model, which we avoid in order to keep the following discussion compact. 

\begin{figure}[tb]
 \centering
 \includegraphics[width=\columnwidth]{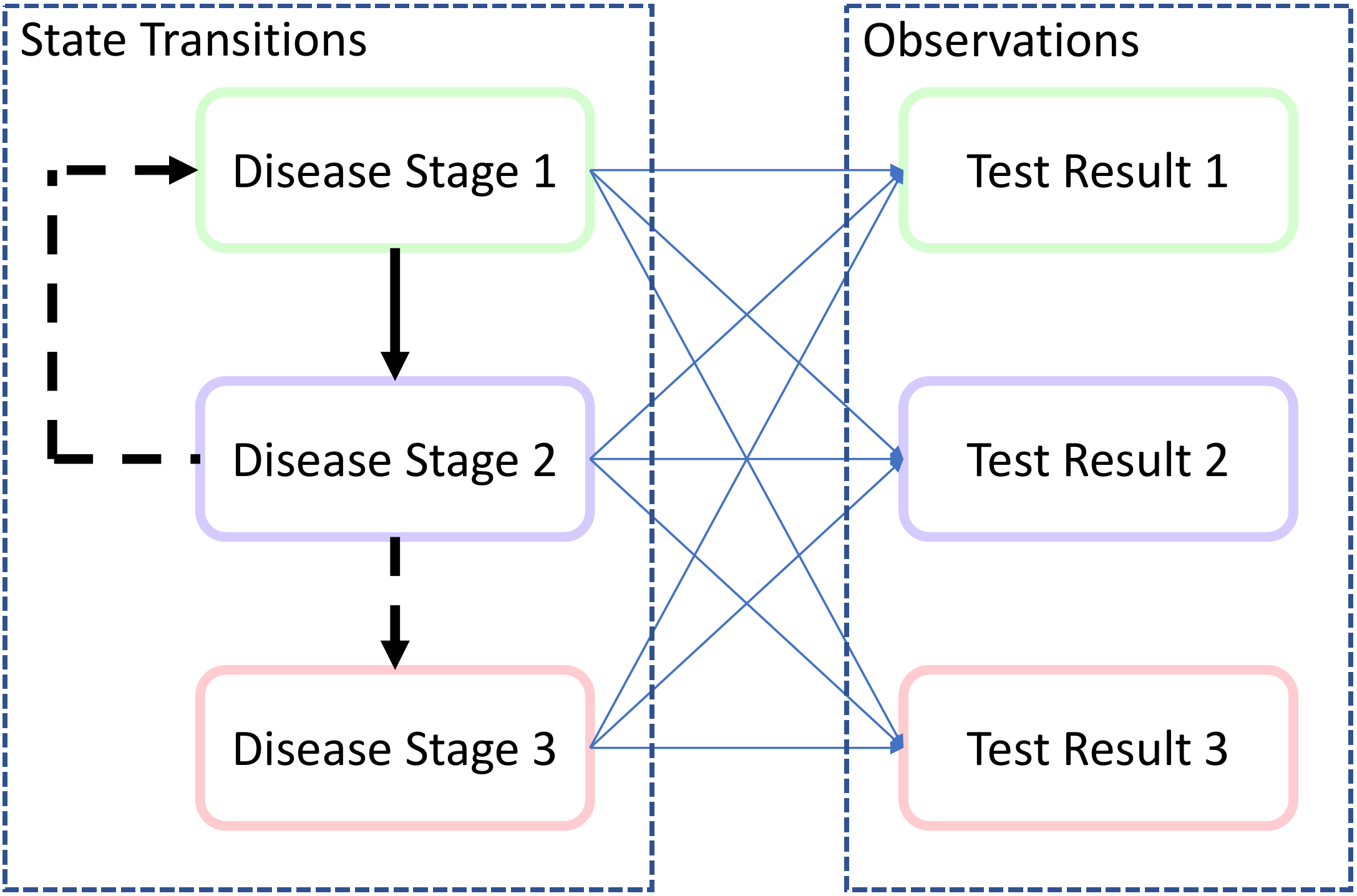}
 \caption{Feasible state transitions and possible test results in healthcare example. Solid arrows for feasible state transitions and observations. Dashed arrows for transitions conditional on treatment and diagnosis decisions.}
 \label{fig:transitions}
\end{figure}

The example, introduced in~\cite{sehr2017tractable}, is set up as follows. The disease encompasses three stages with severity increasing from Stage 1 through Stage 2 to Stage 3, transitions between which are governed by a controlled Markov chain, where $P(a)$ is the transition probability matrix with values $p^a_{ij}$ at row $i$ and column $j$ and $R(a)$ is the observation matrix with elements $r^a_{j\theta}$. All transition and observation probability matrices below are defined similarly. Once our patient enters Stage 3, Stages 1 and 2 are inaccessible for all future times. However, Stage 3 can only be entered through Stage 2, a transition from which to Stage 1 is possible only under costly treatment. The same treatment inhibits transitions from Stage 2 to Stage 3. We have access to the patient state only through imprecise tests, which will result in one of three possible values, each of which is representative of one of the three disease stages. However, these tests are imperfect, with non-zero probability of returning an incorrect disease stage. All possible state transitions and observations are illustrated in Figure~\ref{fig:transitions}.

At each point in time, the current information state $\pi_t$ is available to make one of four possible decisions/actions:
\begin{enumerate}
\item Skip next appointment slot.
\item Schedule new appointment.
\item Order rapid diagnostic test.
\item Apply available treatment.
\end{enumerate}
Skipping an appointment slot results in the patient progressing through the Markov chain describing the transition probabilities of the disease without medical intervention, without new information being available after the current decision epoch. Scheduling an appointment does not alter the patient transition probabilities but provides a low-quality assessment of the current disease stage, which is used to refine the next information state. The third option, ordering a rapid diagnostic test, allows for a high-quality assessment of the patient's state, leading to a more reliable refinement of the next information state than otherwise possible when choosing the previous decision option. The results from this diagnostic test are considered available sufficiently fast so that the patient state remains unchanged under this decision. The remaining option entails medical intervention, allowing probabilistic transition from Stage 2 to Stage 1 while preventing transition from Stage 2 to Stage 3. Transition probabilities $P(a)$, observation probabilities $R(a)$, and stage cost vectors $c(a)$ for each decision are summarized in Table~\ref{tab:numbers}. Additionally, we impose the terminal cost
\begin{align*}
c_N = \begin{bmatrix} 0 & 8 & 60 \end{bmatrix}^T.
\end{align*}

\begin{table*}[htb]
  \caption{Problem data for healthcare decision making example.}
  \centering
  \begin{tabular}{l|ccc}	
    Decision $a$ & Transition Probabilities $P(a)$ & Observation Probabilities $R(a)$ & Cost $c(a)$ \\
    \midrule
\hspace{0.2cm}1: Skip next appointment slot & 
$\begin{bmatrix} 0.80&0.20&0.00\\0.00&0.90&0.10\\0.00&0.00&1.00 \end{bmatrix}$ & 
$\begin{bmatrix} 1/3&1/3&1/3\\1/3&1/3&1/3\\1/3&1/3&1/3 \end{bmatrix}$ & 
$\begin{bmatrix} 0\\5\\5 \end{bmatrix}$ \\[0.5cm]
\hspace{0.2cm}2: Schedule new appointment & 
$\begin{bmatrix} 0.80&0.20&0.00\\0.00&0.90&0.10\\0.00&0.00&1.00 \end{bmatrix}$ & 
$\begin{bmatrix} 0.40&0.30&0.30\\0.30&0.40&0.30\\0.30&0.30&0.40 \end{bmatrix}$ & 
$\begin{bmatrix} 1\\1\\1 \end{bmatrix}$ \\[0.5cm]
\hspace{0.2cm}3: Order rapid diagnostic test & 
$\begin{bmatrix} 1.00&0.00&0.00\\0.00&1.00&0.00\\0.00&0.00&1.00 \end{bmatrix}$ & 
$\begin{bmatrix} 0.90&0.05&0.05\\0.05&0.90&0.05\\0.05&0.05&0.90 \end{bmatrix}$ & 
$\begin{bmatrix} 4\\3\\4 \end{bmatrix}$ \\[0.5cm]
\hspace{0.2cm}4: Apply available treatment & 
$\begin{bmatrix} 0.80&0.20&0.00\\0.75&0.25&0.00\\0.00&0.00&1.00 \end{bmatrix}$ & 
$\begin{bmatrix} 0.40&0.30&0.30\\0.30&0.40&0.30\\0.30&0.30&0.40 \end{bmatrix}$ & 
$\begin{bmatrix} 4\\2\\4 \end{bmatrix}$
  \end{tabular}
  \label{tab:numbers}
\end{table*}

In the solution for the optimal feedback control, the selection of a diagnostic test comes at a cost to the objective criterion and, evidently, serves to refine the information state of the system/patient. It does so without effect on the regulation of the patient other than to improve the information state. Clearly, testing to resolve the state of the patient is part of an optimal strategy in this stochastic setting; but it does take resources. A certainty-equivalent feedback control would assign treatment on the supposition that the patient's state is precisely known. Such a controller would never order a test. The decision to apply a test in the following numerical solution is evidence of duality in receding-horizon stochastic optimal control, viz. SMPC.

\subsection{Computational Results}
\begin{figure*}[h]
  \centering
  \includegraphics[width=\linewidth]{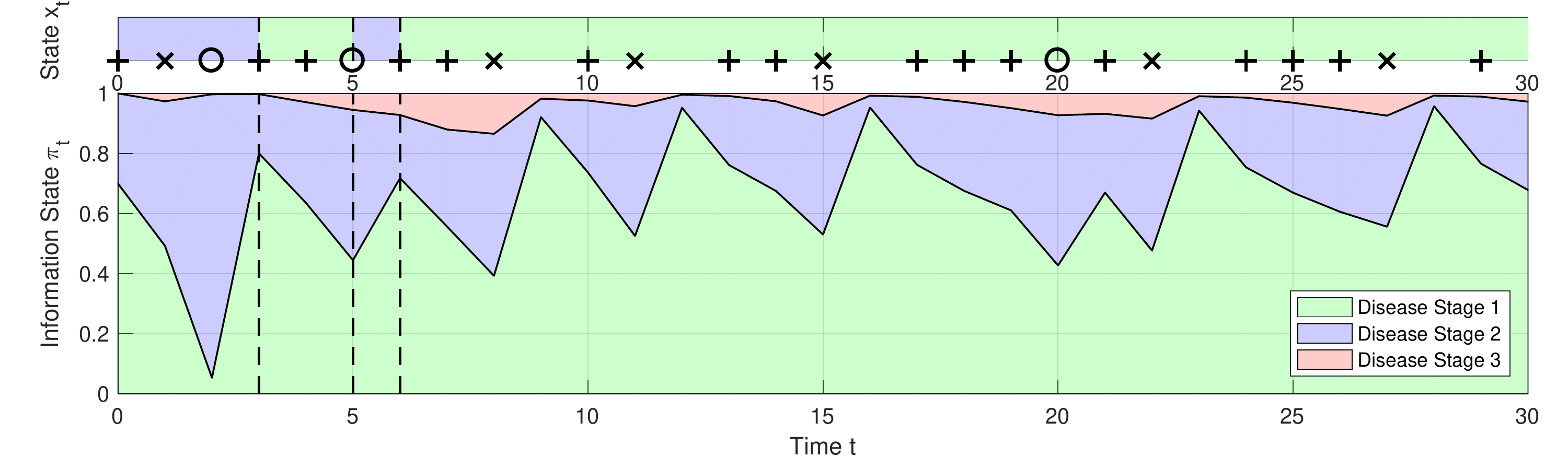}
  \caption{Simulation results for SMPC with horizon $N=4$ and discount factor $\alpha = 0.85$. Top plot displays patient state and transitions, with optimal SMPC decisions based on current information state: appointment (pluses); diagnostic test (crosses); treatment (circles). Bottom plot shows information state evolution. Dashed vertical lines mark time instances of state transitions.}
  \label{fig:simN4}
\end{figure*}

The trade-off between the two principal decision categories -- testing versus treatment, probing versus regulating, exploration versus excitation -- is precisely what is encompassed by duality, which we can include in an optimal sense by solving~(\ref{eq:DP1pomdp}-\ref{eq:DP2pomdp}) and applying the resulting initial policy in receding horizon fashion. This is demonstrated in Figure~\ref{fig:simN4}, which shows simulation results for SMPC with control horizon $N = 4$ and discount factor $\alpha = 0.85$. As anticipated, the stochastic optimal receding horizon policy shows a structure not drastically different from the decision structure motivated above. In particular, diagnostic tests are used effectively to decide on medical intervention. 

In order to apply Theorem~\ref{thm:pomdp} to this particular example, we choose the policy $\tilde{g}(\cdot)$ in Assumption~\ref{assm:pomdp} always to apply medical intervention. Using the worst-case scenario for the expectations in~\eqref{eq:pomdpassm}, which entails transition from Stage 1 to Stage 2 under treatment, we can satisfy Assumption~\ref{assm:pomdp} with $\eta \approx 7$. The computed cost in our simulation is $J_{N}(\pi_0,\mathbf{g}^{\star^{N-1}}) \approx 8.5$. Combined with the discount factor $\alpha = 0.85$, we thus have the upper bound
\begin{align*}
J_{\infty}(\pi_0,\mathbf{g}^{\text{MPC}}) &\leq J_{N}^{\star}(\pi_0) + \frac{\alpha}{1-\alpha}\eta \approx
48.2
\end{align*}
via application of Theorem~\ref{thm:pomdp}. Denoting by $e_j$ the row-vector with entry $1$ in element $j$ and zeros elsewhere, the observed (finite-horizon) cost corresponding with Figure~\ref{fig:simN4} is
\begin{align*}
J_{\infty}^{\text{obs}} = \sum_{k = 0}^{29}e_{x_k}c(\mu_0^N(\pi_k)) \approx 9.2 < 48.2.
\end{align*}

\section{Conclusions}\label{sec:disc}
The central contribution of the paper is the presentation of an SMPC algorithm based on SOOFC. This yields a number of theoretical properties of the controlled system, some of which are simply recognized as the stochastic variants of results from deterministic full-state feedback MPC with their attendant assumptions, including for instance Theorem~\ref{thm:srf} for recursive feasibility. Theorem~\ref{thm:1} is the main stability result in establishing the finiteness of the discounted cost of the SMPC-controlled system. Theorem~\ref{thm:2} and Corollary~\ref{cor:0term} deal with consequent convergence of the state in special cases. 

Performance guarantees of SMPC are made in comparison to performance of the infinite-horizon stochastically optimally controlled system and are presented in Theorem~\ref{thm:bounds} and Corollary~\ref{cor:gamma0v2}. These results extend those of \cite{riggs2012mpc}, which pertain to full-state feedback stochastic optimal control and which therefore do not accommodate duality. Other examples of stochastic performance bounds are mostly restricted to linear systems and, while computable, do not relate to the optimal constrained control. While the formal stochastic results are traceable to deterministic predecessors, the divergence from earlier work is also notable. This concentrates on the use of the information state to accommodate measurements and the exploration of control policy functionals stemming from the Stochastic Dynamic Programming Equation. The resulting output feedback control possesses duality and optimality properties which are either artificially imposed in or absent from earlier approaches.

We have further suggested two potential strategies to ameliorate the computational intractability of the Bayesian filter and SDPE, famous for its \textit{curse of dimensionality}. Firstly, one may use the Particle filter implementation of the Bayesian filter, which has many examples of fast execution for small state dimensions, which with a loss of duality can be combined with scenario methods. This approach is discussed in~\cite{sehr2017particle} as an approximation of the algorithm in this paper. Secondly, we point out that our algorithm becomes computationally tractable for the special case of POMDPs, which may be used either to approximate a nonlinear model or to model a given system in the first place. This strategy inherits the dual nature of our SMPC algorithm for general nonlinear systems.

\appendix
\section{Proofs}
\subsection{Theorem~\ref{thm:1}}
Denote by $M_k$ the discounted  $\mathcal P_N$-cost-to-go,
\begin{align*}
M_k &\triangleq \sum_{j=k}^{k+N-1}{\alpha^{j}c(x_j,g_{j-k}^{\star}(\pi_{j}))} + \alpha^{k+N}c_{N}(x_{k+N})\\
&=\alpha^k \left(\sum_{j=0}^{N-1}{\alpha^{j}c(x_{k+j},g_{j}^{\star}(\pi_{k+j}))} + \alpha^{N}c_{N}(x_{k+N})\right),
\end{align*} 
where $g_j^{\star}(\cdot)$, $j=0,\ldots,N-1$, are the optimal feedback policies in Problem $\mathcal{P}_{N}(\cdot)$. Moreover, define $\mathcal{F}_k$ as the $\sigma$-algebra generated by the initial state $x_0$ with density $\pi_{0\mid -1}$ and the i.i.d. noise sequences $w_j$ and $v_{j}$ for $j = 0,\ldots,k+N-1$. Then $M_k$ is $\mathcal{F}_k$-measurable and $M_k\geq 0$ by non-negativity of stage and terminal cost. Then, 
\begin{multline*}
 \mathbb{E}_0 \left[ M_{k+1} \mid \mathcal{F}_k \right] = \\
\alpha^{k+1} \mathbb{E}_0 \Bigg[\sum_{j=0}^{N-1} \alpha^{j}c(x_{j+k+1},g_{j}^{\star}(\pi_{j+k+1})) \\
+\alpha^{N}c_{N}(x_{k+N+1})  \mid \mathcal{F}_k \Bigg] ,
\end{multline*}
and, by optimality of the policies $g^{\star}_j(\cdot)$ in $\mathcal{P}_N(\cdot)$, 
\begin{align*}
\mathbb{E}_0 [ M_{k+1}& \mid \mathcal{F}_k ] \stackrel{\text{a.s.}}{\leq}
\mathbb{E}_0 \left[ M_k - \alpha^k c(x_{k},g_0^{\star}(\pi_{k})) \right. \\
&- \alpha^{k+N}c_N(x_{k+N}) + \alpha^{k+N}c(x_{k+N},\tilde{g}(\pi_{k+N})) \\
&+\left. \alpha^{k+N+1}c_N(f(x_{k+N},\tilde{g}(\pi_{k+N}),w_{k+N})) \mid \mathcal{F}_k \right],
\end{align*}
where $\tilde{g}(\cdot)$ denotes the terminal feedback policy, specified by Assumptions~\ref{assm:2} and~\ref{assm:policy}, and feasibility follows as in the proof of Theorem~\ref{thm:srf}. Given that
\begin{multline*}
M_k - \alpha^k c(x_{k},g_0^{\star}(\pi_{k})) - \alpha^{k+N}c_N(x_{k+N}) \\
+ \alpha^{k+N}c(x_{k+N},\tilde{g}(\pi_{k+N}))
\end{multline*}
is $\mathcal{F}_k$-measurable, we then have
\begin{align*}
\mathbb{E}_0 [ &M_{k+1} \mid \mathcal{F}_k ] \stackrel{\text{a.s.}}{\leq}
M_k - \alpha^k c(x_{k},g_0^{\star}(\pi_{k})) \\
&- \alpha^{k+N}c_N(x_{k+N}) + \alpha^{k+N}c(x_{k+N},\tilde{g}(\pi_{k+N})) \\
&+\alpha^{k+N+1}\mathbb{E}_0 [c_N(f(x_{k+N},\tilde{g}(\pi_{k+N}),w_{k+N})) \mid \mathcal{F}_k ].
\end{align*}
By Assumption~\ref{assm:policy}, this yields
\begin{align}\label{eq:smc1}
 \mathbb{E}_0 \left[ M_{k+1} \mid \mathcal{F}_k \right] \stackrel{\text{a.s.}}{\leq} M_k - \alpha^k c(x_{k},g_0^{\star}(\pi_{k})).
\end{align}
Taking expectations in~\eqref{eq:smc1} further gives
\begin{align*}
\mathbb{E}_0 \left[ M_{k+1} \right] \leq \mathbb{E}_0 \left[ M_k - \alpha^k c(x_{k},g_0^{\star}(\pi_{k}))\right],
\end{align*}
where $\mathbb{E}_0\left[ M_0 \right] < \infty$ via feasibility of $\pi_0$ for $\mathcal{P}(\cdot)$. By positivity of the stage cost, this yields
\begin{align}\label{eq:L1}
\sup_{k\in\mathbb{N}_0} \mathbb{E}_0 \left[ |M_k| \right] < \infty.
\end{align}
Inequalities~\eqref{eq:smc1} and~\eqref{eq:L1} with non-negativity of the stage cost show that $M_k$ is a non-negative $L^1$-supermartingale on its filtration $\mathcal{F}_k$ and thus, by Doob's Martingale Convergence Theorem (see~\cite{doobclassical}), converges almost surely to a finite random variable,
\begin{align}\label{eq:smc2}
M_k \stackrel{\text{a.s.}}{\to} M_{\infty} < \infty,\text{ as } k\to\infty .
\end{align}
Now define $Z_k$ to be the discounted sample  $\mathcal P_N$ cost-to-go plus the achieved MPC cost at time $k$,
\begin{align*}
Z_k \triangleq M_k + \sum_{j=0}^{k-1}\alpha^j c(x_{j},g_0^{\star}(\pi_{j})) \geq 0.
\end{align*}
Then,
\begin{multline*}
\mathbb{E}_0\left[ Z_{k+1} \mid \mathcal{F}_k \right]\stackrel{\text{a.s.}}{\leq}  \\ 
M_k - \alpha^k c(x_{k},g_0^{\star}(\pi_{k})) + \sum_{j=0}^{k}\alpha^j c(x_{j},g_0^{\star}(\pi_{j})) = Z_k.
\end{multline*}
That is, recognizing that $Z_0=M_0$ so that $\mathbb{E}_0[|Z_0|]<\infty$, $Z_k$ also is a non-negative $L^1$-supermartingale and converges almost surely to a finite random variable 
\begin{align*}
Z_k \stackrel{\text{a.s.}}{\to} Z_{\infty} < \infty,\text{ as } k\to\infty .
\end{align*}
However, by definition of $Z_k$ and~\eqref{eq:smc2}, this implies~\eqref{eq:summability}.
\qed

\subsection{Theorem~\ref{thm:2}}
First proceed as in the proof of Theorem~\ref{thm:1}. By Doob's Decomposition Theorem (see~\cite{doob1990stochastic}) on~\eqref{eq:smc2}, there exists a martingale $\mathcal{M}_k$ and a decreasing sequence $\mathcal{A}_k$ such that $M_k = \mathcal{M}_k + \mathcal{A}_k$, where $A_k\to A_{\infty}$ a.s. by~\eqref{eq:smc2}. Using this decomposition, \eqref{eq:smc1} yields
\begin{multline*}
c(x_{k},g_0^{\star}(\pi_{k})) \leq \alpha^kc(x_{k},g_0^{\star}(\pi_{k}))
\stackrel{\text{a.s.}}{\leq} \\
M_k -  \mathbb{E}_0 \left[ M_{k+1} \mid \mathcal{F}_k \right] = 
\mathcal{A}_k -  \mathbb{E}_0 \left[ \mathcal{A}_{k+1} \mid \mathcal{F}_k \right] \stackrel{\text{a.s.}}{\leq} \\
\mathcal{A}_k -  \mathbb{E}_0 \left[ \mathcal{A}_{\infty} \mid \mathcal{F}_k \right].
\end{multline*}
Taking limits as $k\to\infty$ and re-invoking non-negativity of the stage cost then leads to $c(x_{k},g_0^{\star}(\pi_{k}))\stackrel{}{\to} 0$ a.s., which by the detectability condition on the stage cost (Assumption~\ref{assm:detect}) verifies~\eqref{eq:converge}. 
\qed

\subsection{Theorem~\ref{thm:bounds}}
The optimal value function in the SDPE~\eqref{eq:DP1} satisfies $V_{0}(\pi_0) = J_N(\pi_0,\mathbf{g}^{\star^{N-1}})$, so that optimality of policy $g_0^{\star}(\cdot)$ in Problem $\mathcal{P}_N(\pi_0)$ implies 
\begin{align*}
V_{0}(\pi_0) &= 
\mathbb{E}_{0} \left[ c(x_0,g_0^{\star}(\pi_0)) + \alpha V_{1}(T(\pi_0,y_1,g_0^{\star}(\pi_0))) \right] \\ &+
\alpha\mathbb{E}_{0} \left[ V_{0}(T(\pi_0,y_1,g_0^{\star}(\pi_0))) \right] \\  &- \alpha\mathbb{E}_{0} \left[ V_{0}(T(\pi_0,y_1,g_0^{\star}(\pi_0))) \right],
\end{align*}
which by Assumption~\ref{assm:gamma} yields
\begin{multline}\label{eq:Vineq}
(1 - \alpha \gamma) \mathbb{E}_0 \left[c(x_0,g_0^{\star}(\pi_0)) \right] \leq \\
V_0(\pi_0) - \alpha \mathbb{E}_{0} \left[ V_{0}(T(\pi_0,y_1,g_0^{\star}(\pi_0))) \right] + \alpha \eta.
\end{multline}
Now denote by $J_{\infty}^{M}(\pi_0,\mathbf{g}^{MPC})$ the first $M \in \mathbb{N}_1$ terms of the infinite-horizon cost $J_{\infty}(\pi_0,\mathbf{g}^{MPC})$ subject to the SMPC implementation of policy $g_0^{\ast}(\cdot)$. By~\eqref{eq:Vineq}, we have
\begin{multline*}
(1 - \alpha \gamma) J_{\infty}^{M}(\pi_0,\mathbf{g}^{MPC}) = \\
(1 - \alpha \gamma) \mathbb{E}_0 \left[ \sum_{k=0}^{M-1} \alpha^k c(x_k,g_0^{\star}(\pi_k)) \right] \leq \\
\mathbb{E}_0 \left[ V_0(\pi_0) - \alpha V_0(\pi_1) + \alpha\eta + \alpha V_0(\pi_1) - \alpha^2 V_0(\pi_2)+ \right. \\
\left. \alpha^2\eta +\ldots + 
\alpha^{M-1}V_0(\pi_{M-1}) - \alpha^M V_0(\pi_M) + \alpha^{M} \eta \right],
\end{multline*}
such that
\begin{multline*}
(1 - \alpha \gamma)  J_{\infty}^{M}(\pi_0,\mathbf{g}^{MPC}) \leq 
J_N(\pi_0,\mathbf{g}^{\star^{N-1}})  \\
\quad-\alpha^M \mathbb{E}_0 \left[ J_N(\pi_M,\mathbf{g}^{\star^{N-1}}) \right] + 
\left( \alpha + \ldots +\alpha^{M}\right) \eta,
\end{multline*}
which by non-negativity of the stage cost confirms the right-hand inequality in~\eqref{eq:bounds} in the limit as $M\to\infty$. The left-hand inequality follows directly from optimality.
\qed

\subsection{Corollary~\ref{cor:gamma0v2}}
For conditional densities $\pi_1$ of $x_1$ such that $\pi_1\in\mathcal{C}_{1}$, use optimality and subsequently Assumption~\ref{assm:policyw} to conclude
\begin{align*}
V_{0}(&\pi_1) - V_{1}(\pi_1) \\
=\ &\mathbb{E}_{1} \Bigg[ \left(\sum_{k=0}^{N-1}  \alpha^k c(x_{k+1},g_{k}^{\star}(\pi_{k+1})) + \alpha^Nc_N(x_{N+1}) \right)\\
& -\left(\sum_{k=0}^{N-2}  \alpha^k c(x_{k+1},g_{k+1}^{\star}(\pi_{k+1})) + \alpha^{N-1}c_N(x_{N}) \right) \Bigg] \\ \leq \ & 
\mathbb{E}_{1} [ \alpha^{N-1}c(x_N,\tilde{g}(\pi_N)) \\
&+ \alpha^N c_N(f(x_N,\tilde{g}(\pi_N),w_N)) - \alpha^{N-1}c_N(x_N) ] \\ \leq \ &\eta,
\end{align*}
which by~\eqref{eq:cor:stab} implies $V_0(\pi_k) - V_{1}(\pi_k) \leq\eta$ for $k\in\mathbb{N}_1$. However, this means Assumption~\ref{assm:gamma} is satisfied with $\gamma = 0$ and thus completes the proof by Theorem~\ref{thm:bounds}.
\qed

\bibliographystyle{plain} 
\bibliography{Auto}

\begin{thebibliography}{10}

\bibitem{bertsekas1995dynamic}
D.~P. Bertsekas.
\newblock {\em Dynamic programming and optimal control}.
\newblock Athena Scientific, Belmont, MA, 1995.

\bibitem{bertsekas1978stochastic}
D.~P. Bertsekas and S.~E. Shreve.
\newblock {\em Stochastic optimal control: The discrete time case}, volume~23.
\newblock Academic Press, New York, NY, 1978.

\bibitem{chatterjee2015stability}
D.~Chatterjee and J.~Lygeros.
\newblock On stability and performance of stochastic predictive control
  techniques.
\newblock {\em IEEE Transactions on Automatic Control}, 60(2):509--514, 2015.

\bibitem{chen2003bayesian}
Z.~Chen.
\newblock Bayesian filtering: From {K}alman filters to particle filters, and
  beyond.
\newblock {\em Statistics}, 182(1):1--69, 2003.

\bibitem{copp2014nonlinear}
D.~A. Copp and J.~P. Hespanha.
\newblock Nonlinear output-feedback model predictive control with moving
  horizon estimation.
\newblock In {\em 53rd IEEE Conference on Decision and Control}, pages
  3511--3517, Los Angeles, CA, 2014.

\bibitem{doob1990stochastic}
J.~L. Doob.
\newblock {\em Stochastic processes}.
\newblock John Wiley \& Sons, New York, NY, 1953.

\bibitem{doobclassical}
J.~L. Doob.
\newblock {\em Classical Potential Theory and Its Probabilistic Counterpart}.
\newblock Springer, Berlin, 1984.

\bibitem{elliott2008hidden}
R.~J. Elliott, L.~Aggoun, and J.~B. Moore.
\newblock {\em Hidden {Markov} models: estimation and control}, volume~29.
\newblock Springer Science \& Business Media, 2008.

\bibitem{Feldbaum1965}
A.~A. Fel'dbaum.
\newblock {\em Optimal control systems}.
\newblock Academic Press, New York, NY, 1965.

\bibitem{genceli1996new}
H.~Genceli and M.~Nikolaou.
\newblock New approach to constrained predictive control with simultaneous
  model identification.
\newblock {\em AIChE journal}, 42(10):2857--2868, 1996.

\bibitem{goodwin2014robust}
G.~C. Goodwin, H.~Kong, G.~Mirzaeva, and M.~M. Seron.
\newblock Robust model predictive control: reflections and opportunities.
\newblock {\em Journal of Control and Decision}, 1(2):115--148, 2014.

\bibitem{grune2011nonlinear}
L.~Gr{\"u}ne and J.~Pannek.
\newblock {\em Nonlinear model predictive control}.
\newblock Springer, London, 2011.

\bibitem{grune2008infinite}
L.~Gr{\"u}ne and A.~Rantzer.
\newblock On the infinite horizon performance of receding horizon controllers.
\newblock {\em IEEE Transactions on Automatic Control}, 53(9):2100--2111, 2008.

\bibitem{hernandez1990error}
O.~Hern{\'a}ndez-Lerma and J.~B. Lasserre.
\newblock Error bounds for rolling horizon policies in discrete-time {Markov}
  control processes.
\newblock {\em IEEE Transactions on Automatic Control}, 35(10):1118--1124,
  1990.

\bibitem{kouvaritakis2015stochastic}
B.~Kouvaritakis and M.~Cannon.
\newblock Stochastic model predictive control.
\newblock In John Baillieul and Tariq Samad, editors, {\em Encyclopedia of
  Systems and Control}, pages 1350--1357. Springer, London, 2015.

\bibitem{kouvaritakis2016model}
B.~Kouvaritakis and M.~Cannon.
\newblock {\em Model Predictive Control}.
\newblock Springer, Switzerland, 2016.

\bibitem{BKKUM1986}
P.~R. Kumar and P.~Varaiya.
\newblock {\em Stochastic Systems: Estimation, Identification, and Adaptive
  Control}.
\newblock Prentice-Hall, Englewood Cliffs, NJ, 1986.

\bibitem{marafioti2014persistently}
G.~Marafioti, R.~R. Bitmead, and M.~Hovd.
\newblock Persistently exciting model predictive control.
\newblock {\em International Journal of Adaptive Control and Signal
  Processing}, 28(6):536--552, 2014.

\bibitem{mayne2014model}
D.~Q. Mayne.
\newblock Model predictive control: Recent developments and future promise.
\newblock {\em Automatica}, 50(12):2967--2986, 2014.

\bibitem{mayne2009robust}
D.~Q. Mayne, S.~V. Rakovi{\'c}, R.~Findeisen, and F.~Allg{\"o}wer.
\newblock Robust output feedback model predictive control of constrained linear
  systems: Time varying case.
\newblock {\em Automatica}, 45(9):2082--2087, 2009.

\bibitem{mesbah2016stochastic}
A.~Mesbah.
\newblock Stochastic model predictive control: An overview and perspectives for
  future research.
\newblock {\em IEEE Control Systems Magazine, Accepted}, 2016.

\bibitem{ponomarev1987submersions}
S.~P. Ponomarev.
\newblock Submersions and preimages of sets of measure zero.
\newblock {\em Siberian Mathematical Journal}, 28(1):153--163, 1987.

\bibitem{puterman2014markov}
M.~L. Puterman.
\newblock {\em {Markov} decision processes: discrete stochastic dynamic
  programming}.
\newblock John Wiley \& Sons, 2014.

\bibitem{riggs2012mpc}
D.~J. Riggs and R.~R. Bitmead.
\newblock {MPC} under the hood/sous le capot/unter der haube.
\newblock In {\em 4th IFAC Nonlinear Model Predictive Control Conference},
  pages 363--368, 2012.

\bibitem{schwarm1999chance}
A.~T. Schwarm and M.~Nikolaou.
\newblock Chance-constrained model predictive control.
\newblock {\em AIChE Journal}, 45(8):1743--1752, 1999.

\bibitem{sehr2016sumptus}
M.~A. Sehr and R.~R. Bitmead.
\newblock Sumptus cohiberi: The cost of constraints in {MPC} with state
  estimates.
\newblock In {\em American Control Conference}, pages 901--906, Boston, MA,
  2016.

\bibitem{sehr2017particle}
M.~A. Sehr and R.~R. Bitmead.
\newblock Particle model predictive control: Tractable stochastic nonlinear
  output-feedback {MPC}.
\newblock arXiv:1612.00505. To appear in proc. IFAC World Congress, Toulouse,
  France, 2017.

\bibitem{sehr2017performance}
M.~A. Sehr and R.~R. Bitmead.
\newblock Performance of model predictive control of pomdps.
\newblock ArXiv:1704.07773. Submitted for publication to Proc. 56th IEEE
  Conference on Decision and Control, 2017.

\bibitem{sehr2017tractable}
M.~A. Sehr and R.~R. Bitmead.
\newblock Tractable dual optimal stochastic model predictive control: An
  example in healthcare.
\newblock ArXiv:1704.07770. Submitted for publication to Proc. 1st IEEE
  Conference on Control Technology and Applications, 2017.

\bibitem{simon2006optimal}
D.~Simon.
\newblock {\em Optimal State Estimation: {K}alman, {H}$_{\infty}$, and
  Nonlinear Approaches}.
\newblock John Wiley \& Sons, New York, NY, 2006.

\bibitem{sui2008robust}
D.~Sui, L.~Feng, and M.~Hovd.
\newblock Robust output feedback model predictive control for linear systems
  via moving horizon estimation.
\newblock In {\em American Control Conference}, pages 453--458, Seattle, WA,
  2008.

\bibitem{sunberg2013information}
Z.~Sunberg, S.~Chakravorty, and R.~S. Erwin.
\newblock Information space receding horizon control.
\newblock {\em IEEE Transactions on Cybernetics}, 43(6):2255--2260, 2013.

\bibitem{yan2005incorporating}
J.~Yan and R.~R. Bitmead.
\newblock Incorporating state estimation into model predictive control and its
  application to network traffic control.
\newblock {\em Automatica}, 41(4):595--604, 2005.

\end{thebibliography}

\end{document}